\documentclass[11pt]{article} 

\usepackage[utf8]{inputenc} 

\usepackage{geometry} 
\usepackage{graphicx} 

\usepackage{booktabs} 
\usepackage{array} 
\usepackage{paralist} 
\usepackage{verbatim} 
\usepackage{subfig} 

\usepackage{fancyhdr} 
\usepackage{sectsty}
\allsectionsfont{\sffamily\mdseries\upshape} 

\usepackage[nottoc,notlof,notlot]{tocbibind} 
\usepackage[titles,subfigure]{tocloft} 

\usepackage{amsmath, amssymb, bm}
\usepackage[T1]{fontenc}
\usepackage{textcomp}

\newtheorem{theorem}{Theorem}[section]
\newtheorem{lemma}[theorem]{Lemma}
\newtheorem{proposition}[theorem]{Proposition}

\newtheorem{corollary}[theorem]{Corollary}

\numberwithin{equation}{section} 

\title{Time Evolution of Averaged Limit Shapes of \\ Random Multiple Young Diagrams
\thanks{
This work was supported by JSPS KAKENHI Grant Number JP22K03346} 
}
\author{{Akihito HORA}
\thanks{Department of Mathematics, Faculty of Science, Hokkaido University, Sapporo 
060-0810, Japan; hora@math.sci.hokudai.ac.jp}}

\begin{document}
\maketitle

\begin{abstract}
The branching rule for the tower of wreath products of a finite group by the symmetric groups 
induces a stochastic process on the set of multiple Young diagrams through random 
transitions of boxes of the diagrams between one another. 
We observe dynamical multiple averaged limit shapes resulting from appropriate 
scaling limits, either diffusive or non-diffusive. 
We describe time evolution of the macroscopic multiple averaged limit shapes in terms of 
Voiculescu's $R$-transforms and free L\'evy measures of corresponding Kerov 
transition measures. 
Our microscopic dynamics admits non-exponential pausing time distributions.
\end{abstract}

\section{Introduction}

\begin{figure}[hbt]
\centering
{\unitlength 0.1in%
\begin{picture}(51.2000,6.2700)(4.0000,-8.0000)%
%
\special{pn 8}%
\special{pa 400 400}%
\special{pa 800 800}%
\special{fp}%
\special{pa 800 800}%
\special{pa 1200 400}%
\special{fp}%
\special{pa 440 440}%
\special{pa 480 400}%
\special{fp}%
\special{pa 480 400}%
\special{pa 560 480}%
\special{fp}%
\special{pa 560 480}%
\special{pa 600 440}%
\special{fp}%
\special{pa 600 440}%
\special{pa 680 520}%
\special{fp}%
\special{pa 680 520}%
\special{pa 760 440}%
\special{fp}%
\special{pa 760 440}%
\special{pa 880 560}%
\special{fp}%
\special{pa 880 560}%
\special{pa 920 520}%
\special{fp}%
\special{pa 920 520}%
\special{pa 960 560}%
\special{fp}%
\special{pa 960 560}%
\special{pa 1040 480}%
\special{fp}%
\special{pa 1040 480}%
\special{pa 1080 520}%
\special{fp}%
%
\special{pn 8}%
\special{pa 1200 400}%
\special{pa 1600 800}%
\special{fp}%
\special{pa 1600 800}%
\special{pa 2000 400}%
\special{fp}%
\special{pa 1400 600}%
\special{pa 1440 560}%
\special{fp}%
\special{pa 1440 560}%
\special{pa 1480 600}%
\special{fp}%
\special{pa 1480 600}%
\special{pa 1560 520}%
\special{fp}%
\special{pa 1560 520}%
\special{pa 1640 600}%
\special{fp}%
\special{pa 1640 600}%
\special{pa 1680 560}%
\special{fp}%
\special{pa 1680 560}%
\special{pa 1760 640}%
\special{fp}%
%
\special{pn 8}%
\special{pa 2000 400}%
\special{pa 2400 800}%
\special{fp}%
\special{pa 2400 800}%
\special{pa 2800 400}%
\special{fp}%
\special{pa 2120 520}%
\special{pa 2200 440}%
\special{fp}%
\special{pa 2200 440}%
\special{pa 2240 480}%
\special{fp}%
\special{pa 2240 480}%
\special{pa 2280 440}%
\special{fp}%
\special{pa 2280 440}%
\special{pa 2320 480}%
\special{fp}%
\special{pa 2320 480}%
\special{pa 2400 400}%
\special{fp}%
\special{pa 2400 400}%
\special{pa 2480 480}%
\special{fp}%
\special{pa 2480 480}%
\special{pa 2520 440}%
\special{fp}%
\special{pa 2520 440}%
\special{pa 2560 480}%
\special{fp}%
\special{pa 2560 480}%
\special{pa 2600 440}%
\special{fp}%
\special{pa 2600 440}%
\special{pa 2680 520}%
\special{fp}%
%
\special{pn 8}%
\special{pa 2800 560}%
\special{pa 3120 560}%
\special{fp}%
\special{sh 1}%
\special{pa 3120 560}%
\special{pa 3053 540}%
\special{pa 3067 560}%
\special{pa 3053 580}%
\special{pa 3120 560}%
\special{fp}%
%
\special{pn 20}%
\special{pa 3152 432}%
\special{pa 3184 435}%
\special{pa 3216 440}%
\special{pa 3247 448}%
\special{pa 3277 459}%
\special{pa 3307 471}%
\special{pa 3337 480}%
\special{pa 3368 484}%
\special{pa 3400 483}%
\special{pa 3433 481}%
\special{pa 3465 479}%
\special{pa 3497 481}%
\special{pa 3528 487}%
\special{pa 3559 497}%
\special{pa 3589 510}%
\special{pa 3619 524}%
\special{pa 3650 534}%
\special{pa 3680 535}%
\special{pa 3711 528}%
\special{pa 3741 516}%
\special{pa 3772 506}%
\special{pa 3803 499}%
\special{pa 3816 496}%
\special{fp}%
%
\special{pn 20}%
\special{pa 4112 576}%
\special{pa 4144 575}%
\special{pa 4176 573}%
\special{pa 4208 568}%
\special{pa 4238 559}%
\special{pa 4268 549}%
\special{pa 4299 539}%
\special{pa 4331 533}%
\special{pa 4363 537}%
\special{pa 4389 555}%
\special{pa 4414 575}%
\special{pa 4444 587}%
\special{pa 4473 601}%
\special{pa 4496 616}%
\special{fp}%
%
\special{pn 20}%
\special{pa 4832 496}%
\special{pa 4894 480}%
\special{pa 4924 470}%
\special{pa 4955 460}%
\special{pa 4986 449}%
\special{pa 5016 439}%
\special{pa 5047 433}%
\special{pa 5079 431}%
\special{pa 5111 432}%
\special{pa 5143 431}%
\special{pa 5176 429}%
\special{pa 5207 432}%
\special{pa 5236 445}%
\special{pa 5264 461}%
\special{pa 5295 469}%
\special{pa 5327 472}%
\special{pa 5359 474}%
\special{pa 5390 480}%
\special{pa 5416 488}%
\special{fp}%
%
\special{pn 8}%
\special{pa 3120 400}%
\special{pa 3520 800}%
\special{fp}%
\special{pa 3520 800}%
\special{pa 3920 400}%
\special{fp}%
%
\special{pn 8}%
\special{pa 3920 400}%
\special{pa 4320 800}%
\special{fp}%
\special{pa 4320 800}%
\special{pa 4720 400}%
\special{fp}%
%
\special{pn 8}%
\special{pa 4720 400}%
\special{pa 5120 800}%
\special{fp}%
\special{pa 5120 800}%
\special{pa 5520 400}%
\special{fp}%
%
\special{pn 4}%
\special{pa 3576 520}%
\special{pa 3408 688}%
\special{fp}%
\special{pa 3608 536}%
\special{pa 3432 712}%
\special{fp}%
\special{pa 3648 544}%
\special{pa 3456 736}%
\special{fp}%
\special{pa 3704 536}%
\special{pa 3480 760}%
\special{fp}%
\special{pa 3776 512}%
\special{pa 3504 784}%
\special{fp}%
\special{pa 3544 504}%
\special{pa 3384 664}%
\special{fp}%
\special{pa 3504 496}%
\special{pa 3360 640}%
\special{fp}%
\special{pa 3464 488}%
\special{pa 3336 616}%
\special{fp}%
\special{pa 3416 488}%
\special{pa 3312 592}%
\special{fp}%
\special{pa 3368 488}%
\special{pa 3288 568}%
\special{fp}%
\special{pa 3320 488}%
\special{pa 3264 544}%
\special{fp}%
\special{pa 3288 472}%
\special{pa 3240 520}%
\special{fp}%
\special{pa 3256 456}%
\special{pa 3216 496}%
\special{fp}%
\special{pa 3216 448}%
\special{pa 3192 472}%
\special{fp}%
%
\special{pn 4}%
\special{pa 4416 592}%
\special{pa 4264 744}%
\special{fp}%
\special{pa 4456 600}%
\special{pa 4288 768}%
\special{fp}%
\special{pa 4488 616}%
\special{pa 4312 792}%
\special{fp}%
\special{pa 4392 568}%
\special{pa 4240 720}%
\special{fp}%
\special{pa 4360 552}%
\special{pa 4216 696}%
\special{fp}%
\special{pa 4320 544}%
\special{pa 4192 672}%
\special{fp}%
\special{pa 4248 568}%
\special{pa 4168 648}%
\special{fp}%
\special{pa 4184 584}%
\special{pa 4144 624}%
\special{fp}%
%
\special{pn 4}%
\special{pa 5256 472}%
\special{pa 5024 704}%
\special{fp}%
\special{pa 5296 480}%
\special{pa 5048 728}%
\special{fp}%
\special{pa 5344 480}%
\special{pa 5072 752}%
\special{fp}%
\special{pa 5384 488}%
\special{pa 5096 776}%
\special{fp}%
\special{pa 5224 456}%
\special{pa 5000 680}%
\special{fp}%
\special{pa 5192 440}%
\special{pa 4976 656}%
\special{fp}%
\special{pa 5144 440}%
\special{pa 4952 632}%
\special{fp}%
\special{pa 5096 440}%
\special{pa 4928 608}%
\special{fp}%
\special{pa 5048 440}%
\special{pa 4904 584}%
\special{fp}%
\special{pa 4976 464}%
\special{pa 4880 560}%
\special{fp}%
\special{pa 4904 488}%
\special{pa 4856 536}%
\special{fp}%
%
\special{pn 8}%
\special{pa 720 400}%
\special{pa 720 240}%
\special{fp}%
\special{sh 1}%
\special{pa 720 240}%
\special{pa 700 307}%
\special{pa 720 293}%
\special{pa 740 307}%
\special{pa 720 240}%
\special{fp}%
\special{pa 880 240}%
\special{pa 880 400}%
\special{fp}%
\special{sh 1}%
\special{pa 880 400}%
\special{pa 900 333}%
\special{pa 880 347}%
\special{pa 860 333}%
\special{pa 880 400}%
\special{fp}%
%
\special{pn 8}%
\special{pa 1520 400}%
\special{pa 1520 240}%
\special{fp}%
\special{sh 1}%
\special{pa 1520 240}%
\special{pa 1500 307}%
\special{pa 1520 293}%
\special{pa 1540 307}%
\special{pa 1520 240}%
\special{fp}%
\special{pa 1680 240}%
\special{pa 1680 400}%
\special{fp}%
\special{sh 1}%
\special{pa 1680 400}%
\special{pa 1700 333}%
\special{pa 1680 347}%
\special{pa 1660 333}%
\special{pa 1680 400}%
\special{fp}%
%
\special{pn 8}%
\special{pa 2320 400}%
\special{pa 2320 240}%
\special{fp}%
\special{sh 1}%
\special{pa 2320 240}%
\special{pa 2300 307}%
\special{pa 2320 293}%
\special{pa 2340 307}%
\special{pa 2320 240}%
\special{fp}%
\special{pa 2480 240}%
\special{pa 2480 400}%
\special{fp}%
\special{sh 1}%
\special{pa 2480 400}%
\special{pa 2500 333}%
\special{pa 2480 347}%
\special{pa 2460 333}%
\special{pa 2480 400}%
\special{fp}%
\end{picture}}%
\vspace{-5mm}
\caption{Random multi-diagrams to multi-interfaces: the case of $|\widehat{T}| =3$}
\label{fig:1-1}
\end{figure}
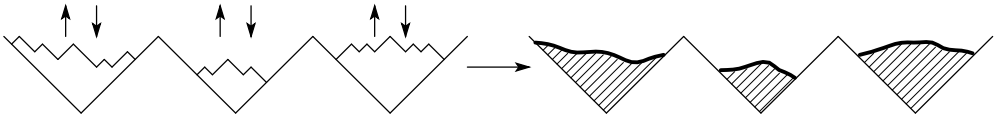

In the present note, we observe a certain time evolution of 
multi-interfaces as indicated in Figure~\ref{fig:1-1} (right). 
Sufficiently zoomed in, each interface looks like a zigzag profile of a Young diagram. 
Moreover, the profiles fluctuate by exchanging boxes at random between one another. 
The random structure of our model comes from the branching rule for a tower of 
wreath product groups. 
Actually, a broad aim of our project is to study asymptotic behavior of the brancing 
rule through probability models. 
Let us settle the model precisely. 

We consider the wreath product of a finite group $T$ by the symmetric group 
$\mathfrak{S}_n$ of degree $n$: 
\[
\mathfrak{S}_n(T) = T \wr \mathfrak{S}_n = T^n \rtimes \mathfrak{S}_n. 
\]
The set of conjugacy classes and the set of equivalence classes of the irreducible 
representations of $T$ are denoted by $[T]$ and $\widehat{T}$ respectively. 
The symbol $\mathbb{Y}_n$ denotes the set of Young diagrams with $n$ boxes. 
Set $\mathbb{Y} = \bigsqcup_{n=0}^\infty \mathbb{Y}_n$ where $\mathbb{Y}_0$ 
consists of empty diagram $\varnothing$. 
The set of conjugacy classes $[\mathfrak{S}_n(T)]$ and the set of equivalence classes 
of the irreducible representations $(\mathfrak{S}_n(T))^{\wedge}$ of 
$\mathfrak{S}_n(T)$ are parametrized by the $|[T]| = |\widehat{T}|$-tuples of 
Young diagrams with the numbers of boxes summing up to $n$, denoted by 
$\mathbb{Y}_n([T])$ and $\mathbb{Y}_n(\widehat{T})$ respectively: 
\begin{align*}
&[\mathfrak{S}_n(T)] \cong 
\mathbb{Y}_n([T]) = \bigl\{ \rho = (\rho_\theta)_{\theta\in [T]} \,\big|\, 
\rho_\theta\in\mathbb{Y}, \; \sum_{\theta\in [T]} |\rho_\theta| =n \bigr\}, \\ 
&(\mathfrak{S}_n(T))^{\wedge} \cong 
\mathbb{Y}_n(\widehat{T}) = \bigl\{ \lambda = (\lambda^\zeta)_{\zeta\in \widehat{T}} 
\,\big|\, \lambda^\zeta\in\mathbb{Y}, \; \sum_{\zeta\in \widehat{T}} 
|\lambda^\zeta| =n \bigr\} 
\end{align*}
where $|\cdot|$ denotes the number of boxes of a Young diagram. 
For $k, n\in\mathbb{N}$ such that $k<n$, the canonical inclusion 
$\iota_{n,k} : \mathfrak{S}_k(T) \longrightarrow \mathfrak{S}_n(T)$ 
\begin{align*}
&x\in T^k \longmapsto (x, e_T, \cdots, e_T) \in T^k\times T^{n-k} = T^n,  \\ 
&\sigma\in\mathfrak{S}_k \longmapsto \sigma (k+1)\cdots (n) \in \mathfrak{S}_n 
\end{align*}
induces the following inclusion, which is denoted again by 
$\iota_{n,k} : \mathbb{Y}_k([T]) \longrightarrow \mathbb{Y}_n([T])$. 
For $\rho = (\rho_\theta)_{\theta\in [T]} \in \mathbb{Y}_k([T])$, 
\[
\iota_{n,k}\rho = (\Tilde{\rho}_\theta)_{\theta\in [T]} \in \mathbb{Y}_n([T]) 
\qquad \text{where} \quad \Tilde{\rho}_\theta = 
\begin{cases} \rho_{\{e_T\}\sqcup (1^{n-k})}, & \theta = \{e_T\} \\ 
\rho_\theta, & \text{otherwise}. \end{cases} 
\] 
For simplicity we often write as $\mathfrak{S}_k(T)$ instead of $\iota_{n,k}\mathfrak{S}_k(T)$ 
to have $\mathfrak{S}_k(T) \subset \mathfrak{S}_n(T)$. 

For $\lambda = (\lambda^\zeta)_{\zeta\in\widehat{T}} \in \mathbb{Y}_n(\widehat{T})$ 
and $\nu = (\nu^\zeta)_{\zeta\in\widehat{T}} \in \mathbb{Y}_{n-1}(\widehat{T})$, 
we write as $\nu\nearrow\lambda$ or $\lambda\searrow\nu$ when $\lambda$ is obtained 
by adding a box at one entry of $\nu$. 
The uniquely determined entry label is then denoted by $\zeta_{\nu\lambda}\in\widehat{T}$. 
Let $\pi^\lambda$ be an irreducible unitary representation of $\mathfrak{S}_n(T)$ 
corresponding to $\lambda\in\mathbb{Y}_n(\widehat{T})$. 
The dimension of $\pi^\lambda$ is given by 
\begin{equation}\label{eq:1-7}
\dim\lambda = \frac{n!}{\prod_{\zeta\in\widehat{T}} |\lambda^\zeta|!} 
\prod_{\zeta\in\widehat{T}} (\dim \zeta)^{|\lambda^\zeta|} \dim\lambda^\zeta.
\end{equation}
We have irreducible decomposition with dimension counting: 
\begin{equation}\label{eq:1-8}
\mathrm{Res}^{\mathfrak{S}_n(T)}_{\mathfrak{S}_{n-1}(T)} \pi^\lambda 
\cong \!\!\bigoplus_{\nu\in\mathbb{Y}_{n-1}(\widehat{T}) : \, \nu\nearrow\lambda} 
[\dim\zeta_{\nu\lambda}] \pi^\nu, \quad 
\dim\lambda = 
\!\!\sum_{\nu\in\mathbb{Y}_{n-1}(\widehat{T}) : \, \nu\nearrow\lambda} 
\dim\zeta_{\nu\lambda} \dim\nu
\end{equation}
and reciprocally, 
\begin{equation}\label{eq:1-9}
\mathrm{Ind}^{\mathfrak{S}_n(T)}_{\mathfrak{S}_{n-1}(T)} \pi^\nu 
\cong \!\!\bigoplus_{\mu\in\mathbb{Y}_n(\widehat{T}) : \, \nu\nearrow\mu} 
[\dim\zeta_{\nu\mu}] \pi^\mu, \quad 
n |T| \dim\nu = 
\!\!\sum_{\mu\in\mathbb{Y}_n(\widehat{T}) : \, \nu\nearrow\mu} 
\dim\zeta_{\nu\mu} \dim\mu.
\end{equation}
\eqref{eq:1-8} and \eqref{eq:1-9} produce two stochastic matrices by 
\[
P^\downarrow_{\lambda\nu} = \frac{\dim\zeta_{\nu\lambda}\dim\nu}{\dim\lambda}, 
\qquad 
P^\uparrow_{\nu\mu} = \frac{\dim\zeta_{\nu\mu}\dim\mu}{n |T| \dim\nu}
\] 
and thereby a stochastic matrix of degree $|\mathbb{Y}_n(\widehat{T})|$ 
\begin{equation}\label{eq:1-11}
P = P^{(n)} = P^\downarrow P^\uparrow.
\end{equation}
More explicitly, the entries of $P$ of \eqref{eq:1-11} are given by 
\[
P_{\lambda\mu} = \begin{cases} \dfrac{1}{n}\, \dfrac{|\lambda^\zeta|}{|\mu^\eta|}
\,\dfrac{\dim\mu^\zeta \dim\mu^\eta}{\dim\lambda^\zeta \dim\lambda^\eta} 
\,\dfrac{(\dim\eta)^2}{|T|}, & \exists \zeta, \eta \text{ s.t. } \lambda^\zeta\searrow\mu^\zeta, 
\ \lambda^\eta\nearrow\mu^\eta \\ 
\dfrac{1}{n}\, \dfrac{\dim\mu^\zeta}{\dim\lambda^\zeta}\, \dfrac{(\dim\zeta)^2}{|T|}, 
& \lambda\neq\mu; \  \forall\zeta, \; |\lambda^\zeta|=|\mu^\zeta| \\ 
\dfrac{1}{n|T|} \sum_\zeta (\dim\zeta)^2\cdot \#\{ \text{peaks of }\lambda^\zeta\}, 
& \lambda =\mu \\ 
\ 0, & \text{otherwise}. \end{cases}
\]
The Markov chain induced by \eqref{eq:1-11}, called Res-Ind chain on 
$\mathbb{Y}_n(\widehat{T})$, is symmetric (i.e. satisfies detailed balance) 
with respect to the Plancherel measure 
\[
M_{\mathrm{Pl}}(\lambda) = \frac{(\dim\lambda)^2}{n! |T|^n}, \qquad 
M_{\mathrm{Pl}}(\lambda) P_{\lambda\mu} = M_{\mathrm{Pl}}(\mu) P_{\mu\lambda}, 
\qquad \lambda, \mu\in \mathbb{Y}_n(\widehat{T}).
\] 
The procedure of defining a Res-Ind chain works well for a general pair $(G, H)$ of 
a finite group and its subgroup. 
Decomposition of $\mathrm{Ind}\circ\mathrm{Res} \,\pi^\lambda$ for wreath products 
was considered in \cite{Oka}. 
Res-Ind chains were introduced by Fulman to apply to the Stein method for algebraic 
central limit theorems in \cite{Ful04}, \cite{Ful05}. 
See \cite{BoOl09} for an application to infinite-dimensional diffusions. 
We have treated Res-Ind chains in a series of works \cite{Hor15}, \cite{Hor16}, \cite{Hor20} 
and \cite{Hor24} preceding the present note. 

Iteration of \eqref{eq:1-8} and \eqref{eq:1-9} along $n$ produces the branching graph 
for the towers of wreath product groups. 
Our series of works \cite{HiHiHo09}, \cite{HoHiHi08} and \cite{HoHi14} treat related 
harmonic analysis where asymptotic analysis is developed in a different scaling regime from 
this note. 
Interesting probability models on such a branching graph are studied in \cite{Str24}. 
A static model on fluctuations for growing multiple Young diagrams is studied in \cite{Sni06}. 

Substituting a counting process for discrete time $k$ of Res-Ind chain 
$(Z_k^{(n)})_{k\in\{0,1,2,\dots\}}$, we obtain a continuous time process on 
$\mathbb{Y}_n(\widehat{T})$. 
For that purpose, let $\psi$ be a pausing (or holding) time distribution, which is a 
probability on $(0, \infty)$. 
Take an IID sequence $(\epsilon_j)_{j\in\mathbb{N}}$, each obeying $\psi$, and set 
\[
N_s = \begin{cases} 0, & s< \epsilon_1 \\ 
j, & \epsilon_1+\cdots+\epsilon_j\leqq s< \epsilon_1+\cdots+\epsilon_{j+1}, \end{cases} 
\qquad N_0 =0 \ \text{a.s.}
\] 
The continuous time process $(X_s^{(n)})_{s\geqq 0}$ obtained by 
\begin{equation}\label{eq:1-14}
X_s^{(n)} = Z_{N_s}^{(n)}, \qquad s\geqq 0 
\end{equation}
describes our microscopic random dynamics. 
The distribution at time $s$ is expressed as 
\begin{equation}\label{eq:1-15}
\mathbb{P}( X_s^{(n)} =\mu \,|\, X_0^{(n)} = \lambda) = 
\sum_{j=0}^\infty (P^{(n)j})_{\lambda\mu} \int_{[0, s]} 
\psi\bigl( (s-u, \infty)\bigr) \psi^{\ast j}(du), \qquad 
\lambda, \mu\in \mathbb{Y}_n(\widehat{T})
\end{equation}
where $\psi^{\ast j}$ is the $j$th convolution of $\psi$. 
See \cite{Hor20}, \cite{Hor24} for \eqref{eq:1-15}
\footnote{The assumption of $\psi((0,\infty)) >0$ in \cite{Hor20} and \cite{Hor24} 
is insufficient and replaced by $\psi(\{0\}) =0$ for $N_0 =0$ a.s.}. 

Transfer from microscopic random multi-diagrams to macroscopic multi-interfaces as 
Figure~\ref{fig:1-1} is formulated through appropriate space--time scaling limits. 
In the spatial direction, $1/\sqrt{n}$ rescale is nice for the limit shape of Young diagrams. 
We set 
\begin{equation}\label{eq:1-16}
(\lambda^\zeta)^{\sqrt{n}} (x) = \frac{1}{\sqrt{n}} \lambda^\zeta (\sqrt{n} x), \qquad 
\lambda = (\lambda^\zeta)_{\zeta\in\widehat{T}} \in\mathbb{Y}_n(\widehat{T}) 
\end{equation}
where $y = \lambda^\zeta (x)$ is the profile of $\lambda^\zeta$ displayed in the $xy$-plane. 
In the temporal direction, letting $t$ be a macroscopic time, we set $s= t \tau_n$ 
where $(\tau_n)$ is a positive sequence such that $\lim_{n\to\infty}\tau_n = \infty$. 
Later the diverging order $\tau_n$ will be taken in several ways. 
Compared with spatial $1/\sqrt{n}$ rescale, the case of $\tau_n =n$ gives a diffusive limit. 
In our project, we investigate the behavior of random multi-diagrams 
\begin{equation}\label{eq:1-17}
\bigl( ( X_{t\tau_n}^{(n)\, \zeta})^{\sqrt{n}} \bigr)_{\zeta\in\widehat{T}} 
\end{equation}
as $n\to\infty$ and time evotution (i.e. concrete $t$-dependence) of limiting objects 
by setting appropriate assumptions on initial ensembles (i.e. distributions of $X_0^{(n)}$) 
and pausing time distributions. 
In this note, we derive the time evolution of some \textit{averaged} quantities for 
\eqref{eq:1-17} by applying harmonic analysis of wreath product groups. 
The results for averaged limit shapes in the case of stable pausing time distributions are new, 
as continuation from \cite{Hor20}, even if they are restricted to trivial $T$. 
A more detailed aspect of dymanical multi-interfaces as a \textit{concentration} phenomenon 
(i.e. law of large numbers) will be discussed in a separate note. 

The subsequent sections are organized as follows. 
The next \S2 is devoted to reviewing some necessary facts including the character formula 
for wreath product groups. 
In \S3, we state evolution of the averaged areas (i.e. numbers of boxes) of rescaled 
multi-diagrams, which are the simplest macroscopic quantities resulting from interaction 
between multi-diagrams. 
In \S4, we investigate $t$-dependence of the transition measures to fully describe the 
averaged dynamical limit shapes. 
Computation of their free cumulants is essential. 
Examples of typical ensembles are constructed from characters of the infinite wreath product 
$\mathfrak{S}_\infty(T)$ in \S5. 

\section{Preliminaries}

\subsection{Character formula for wreath product groups}

Character formula of a wreath product group is well-known. 
We will consider asymptotic behavior in which the group size and irreducible representations 
grow while type of a conjugacy class is fixed, namely in the dual approach due to 
Kerov--Olshanski. 
Keeping that purpose in mind, we recall the presentations used in \cite{HoHi14}.  

Let $k, n\in\mathbb{N}$ such that $k<n$, 
$\rho = (\rho_\theta)_{\theta\in [T]} \in \mathbb{Y}_k([T])$ and 
$\lambda = (\lambda^\zeta)_{\zeta\in\widehat{T}} \in \mathbb{Y}_n(\widehat{T})$. 
The irreducible character value $\chi^\lambda_{\iota_{n,k}\rho}$ is given as follows. 
Let $\mathrm{rows}(\rho)$ denote the set of rows of $\rho$, where each $\rho_\theta$ 
is decomposed into rows and every row is distinguished in (possibly) equal length. 
A map $r : \mathrm{rows}(\rho) \longrightarrow \widehat{T}$ gives a way to assign 
an element of $\widehat{T}$ to each row in $\mathrm{rows}(\rho)$. 
Among these maps, we specify the one satisfying: 
\[
\text{for any } \zeta\in\widehat{T}, \quad (\text{the sum of lengths of rows in } 
r^{-1}(\zeta))\  \leqq |\lambda^\zeta| 
\] 
and calle it admissible with respect to $\lambda$ (temporally). 
The Young diagram consisting of the rows in $r^{-1}(\zeta)$ is denoted again by 
$r^{-1}(\zeta)$, whose size is $|r^{-1}(\zeta)|$. 
Moreover, $r^{-1}(\zeta)\cap\rho_\theta$ is the Young diagram consisting of the rows in 
$\rho_\theta$ to which label $\zeta$ is assigned. 
As usual, $l(\,\cdot\,)$ denotes the number of rows of a Young diagram. 
We have (see Eq.(1.12) in \cite{HoHi14}) 
\begin{multline}\label{eq:2-1-2}
\chi^\lambda_{\iota_{n,k}\rho} = 
\sum_{\text{admissible } r :\, \mathrm{rows}(\rho) \to \widehat{T}} 
\frac{(n-k)!}{\prod_{\zeta\in\widehat{T}} (|\lambda^\zeta|-|r^{-1}(\zeta)|)!} 
\prod_{\zeta\in\widehat{T}} \Bigl\{ (\dim\zeta)^{|\lambda^\zeta|-|r^{-1}(\zeta)|} \\ 
\bigl( \prod_{\theta\in [T]} (\chi^\zeta_\theta)^{l(r^{-1}(\zeta)\cap\rho_\theta)}\bigr) 
\chi^{\lambda^\zeta}_{r^{-1}(\zeta)\sqcup (1^{|\lambda^\zeta|-|r^{-1}(\zeta)|})}\Bigr\}
\end{multline}
where $(\chi^\zeta_\theta)_{\zeta, \theta}$ is the character table of $T$. 
Let us set
\footnote{Notation for falling factorial:\ $n^{\downarrow k} = n(n-1) \cdots (n-k+1)$}: 
for $\tau\in\mathbb{Y}$ and $k\in\mathbb{N}$, 
\begin{equation}\label{eq:2-1-3}
\Sigma_\tau (\nu) = |\nu|^{\downarrow |\tau|} 
\frac{\chi^\nu_{(\tau, 1^{|\nu|-|\tau|})}}{\dim\nu} 
\quad \bigl( =0 \text{ if } |\nu|<|\tau| \bigr), \quad \Sigma_k(\nu) = \Sigma_{(k)}(\nu), 
\qquad \nu\in\mathbb{Y}.
\end{equation}
Using \eqref{eq:1-7} and \eqref{eq:2-1-3}, we rewrite \eqref{eq:2-1-2} as 
\begin{equation}\label{eq:2-1-4}
\frac{\chi^\lambda_{\iota_{n,k}\rho}}{\dim\lambda} = 
\frac{1}{n^{\downarrow k}} \sum_{r :\, \mathrm{rows}(\rho) \to \widehat{T}} 
\prod_{\zeta\in\widehat{T}} \Bigl\{ 
\frac{\Sigma_{r^{-1}(\zeta)}(\lambda^\zeta)}{(\dim\zeta)^{|r^{-1}(\zeta)|}} 
\prod_{\theta\in [T]} (\chi^\zeta_\theta)^{l(r^{-1}(\zeta)\cap\rho_\theta)}\Bigr\}. 
\end{equation}
Note that non-admissible terms vanish by the notation of \eqref{eq:2-1-3}. 
In particular, for $\theta\in [T]$ and $k\in\mathbb{N}$, let $(k)_\theta$ denote the 
$(\theta, k)$-cycle type in $\mathbb{Y}_k([T])$, i.e. the $\theta$-entry is $k$-cycle $(k)$ and 
the other entries are all $\varnothing$. 
Then, $\mathrm{rows}((k)_\theta)$ is a singleton. 
The character formula \eqref{eq:2-1-4} yields for $\lambda\in\mathbb{Y}_n(\widehat{T})$ 
\begin{equation}\label{eq:2-1-5}
\frac{\chi^\lambda_{\iota_{n,k}(k)_\theta}}{\dim\lambda} 
= \frac{1}{n^{\downarrow k}} \sum_{\zeta\in\widehat{T}} 
\frac{\chi^\zeta_\theta}{(\dim\zeta)^k} \Sigma_k(\lambda^\zeta). 
\end{equation}

\subsection{Eigenvalues and eigenvectors of a Res-Ind chain}

In general, the transition matrix \eqref{eq:1-11} of a Res-Ind chain for finite group $G$ and its 
subgroup $H$ has eigenvectors consisting of each column of the normalized character table of $G$. 
To be precise, consider normalized irreducible characters of $G$ 
and their values at an element in conjugacy class $C$ of $G$: 
$\chi^\xi_C / \dim\xi$ ($\xi\in\widehat{G}$). 
We have 
\begin{equation}\label{eq:2-2-1}
P v_C = \frac{|C\cap H|}{|C|} \, v_C \qquad \text{where} \quad 
v_C = \Bigl( \frac{\chi^\xi_C}{\dim\xi} \Bigr)_{\xi\in\widehat{G}}.
\end{equation}
We refer to \cite{Hor15} and \cite{Hor16} for \eqref{eq:2-2-1}. 
Now let $G = \mathfrak{S}_n(T)$, $H=\mathfrak{S}_{n-1}(T)$. 
Conjugacy classes of a wreath product group are described well through standard decomposition 
into basic elements (see e.g. \cite{HiHi05}). 
Let the type of $C$ be $\rho = (\rho_\theta)_{\theta\in [T]}\in\mathbb{Y}_k([T])$. 
Then
\footnote{Notation for a Young diagram by the multiplicities of parts: 
$\nu = (1^{m_1(\nu)} 2^{m_2(\nu)} 3^{m_3(\nu)}\cdots)$.}, 
\begin{equation}\label{eq:2-2-2}
\frac{|C\cap H|}{|C|} = 1- \frac{k- m_1(\rho_{\{e_T\}})}{n}
\end{equation}
holds. 
Applying \eqref{eq:2-2-1} and \eqref{eq:2-2-2} to the transition matrix $P^{(n)}$ of 
\eqref{eq:1-11}, we have 
\begin{equation}\label{eq:2-2-2.5}
P^{(n)} \Bigl( \frac{\chi^\lambda_{\iota_{n,k}\rho}}{\dim\lambda}
\Bigr)_{\lambda\in\mathbb{Y}_n(\widehat{T})} = 
\Bigl( 1- \frac{k- m_1(\rho_{\{e_T\}})}{n}\Bigr) 
\Bigl( \frac{\chi^\lambda_{\iota_{n,k}\rho}}{\dim\lambda}
\Bigr)_{\lambda\in\mathbb{Y}_n(\widehat{T})}, 
\end{equation}
in particular 
\begin{equation}\label{eq:2-2-3}
P^{(n)} \Bigl( \frac{\chi^\lambda_{\iota_{n,k}(k)_\theta}}{\dim\lambda}
\Bigr)_{\lambda\in\mathbb{Y}_n(\widehat{T})} = \bigl( 1- \frac{k}{n}\bigr) 
\Bigl( \frac{\chi^\lambda_{\iota_{n,k}(k)_\theta}}{\dim\lambda}
\Bigr)_{\lambda\in\mathbb{Y}_n(\widehat{T})} 
\end{equation}
for $k\geqq 2$, or for $k=1$ and $\theta \neq \{e_T\}$. 
If $k=1$ and $\theta = \{e_T\}$, \eqref{eq:2-2-2.5} implies 
(by using the column vector with $1$ as all entries) an obvious equality 
\[ 
P^{(n)} (1)_{\lambda\in\mathbb{Y}_n(\widehat{T})} = 
(1)_{\lambda\in\mathbb{Y}_n(\widehat{T})}.
\] 

\subsection{Kerov transition measure, Kerov polynomials, etc.}

In this note, the profile of Young diagram $\nu$ put in the $xy$-plane is denoted by the 
same symbol $y=\nu(x)$. 
As for its size, we have 
\[ 
\int_\mathbb{R} (\nu(x)-|x|)\, dx = 2 |\nu|.
\] 
The coordinates of peaks and valleys of a profile $x_1 <y_1<\cdots < y_{r-1}<x_r$ 
satisfy \ $x_i, y_i\in\mathbb{Z}$ and $\sum_i x_i = \sum_i y_i$. 
If $x_i$ and $y_i$ are not necessarily integers, such an interlacing sequence determines 
a profile of a rectangular diagram. 
More generally, an element of 
\[ 
\mathbb{D} = \{ \omega : \mathbb{R} \longrightarrow \mathbb{R} \,|\, 
|\omega(x)-\omega(y)| \leqq |x-y|, \ \omega(x) =|x| \text{ for sufficiently large } |x|\}
\] 
is a continuous diagram. 
A compactly supported probability $\mathfrak{m}_\omega$ on $\mathbb{R}$ is 
assigned to $\omega\in\mathbb{D}$, which is referred to as the Kerov transition 
measure of $\omega$, by 
\begin{equation}\label{eq:2-3-3}
\frac{1}{z} \exp\Bigl\{ \int_\mathbb{R} \frac{1}{x-z} \Bigl( 
\frac{\omega(x)-|x|}{2}\Bigr)^\prime dx \Bigr\} = 
\int_\mathbb{R} \frac{1}{z-x}\, \mathfrak{m}_\omega(dx), \qquad z\in\mathbb{C}^+.
\end{equation}
A continuous diagram is fully encoded by its transition measure. 
The correspondence of $\omega \leftrightarrow \mathfrak{m}_\omega$ is called 
the Markov transform. 
If $\omega$ is a rectanglular diagram determined by $x_1 <y_1<\cdots < y_{r-1}<x_r$, 
$\mathfrak{m}_\omega$ is an atomic probability supported by $x_i$'s, and \eqref{eq:2-3-3} 
is reduced to 
\[ 
\frac{(z-y_1)\cdots (z-y_{r-1})}{(z-x_1)\cdots (z-x_r)} = 
\sum_{i=1}^r \frac{\mathfrak{m}_\omega(x_i)}{z-x_i}.
\] 
In this note, we do not treat transition measures with noncompact support nor corresponding 
continuous diagrams. 

The algebra of polynomial functions of Young diagrams (Kerov--Olshanski algebra) 
$\mathbb{A}$ contains several sorts of generators. 
When $\nu\in\mathbb{Y}$ has the interlacing coordinates $x_1<y_1\cdots <y_{r-1}<x_r$, 
set $q_k(\nu) = \sum_i x_i^k - \sum_i y_i^k$. 
The system $\{q_k\}_{k\in\mathbb{N}}$ generates $\mathbb{A}$. 
Regarding $q_k$ as a homogeneous element of degree $k$ in $\mathbb{A}$, we define 
weight degrees in $\mathbb{A}$: $\mathrm{wt} q_k =k$. 
The moments and free cumulants of transition measures, 
$\{ M_k(\mathfrak{m}_\nu)\}_{k\in\mathbb{N}}$ and 
$\{ R_k(\mathfrak{m}_\nu)\}_{k\in\mathbb{N}}$ respectively, are generator sets of 
$\mathbb{A}$ also, in which 
$\mathrm{wt} M_k(\mathfrak{m}_\nu) = \mathrm{wt} R_k(\mathfrak{m}_\nu) =k$. 
The two systems are connected by the free cumulant-moment formula. 
$\{\Sigma_\tau\}_{\tau\in\mathbb{Y}}$ is a linear basis of $\mathbb{A}$ 
with $\mathrm{wt}\Sigma_\tau = |\tau|+l(\tau)$. 
In $\mathbb{A}$, $\{\Sigma_k(\nu)\}_{k\in\mathbb{N}}$ and 
$\{ R_k(\mathfrak{m}_\nu)\}_{k\in\mathbb{N}}$ are related by the Kerov polynomials: 
for $k\in\mathbb{N}$ such that $k\geqq 3$, 
\begin{equation}\label{eq:2-3-5}
\Sigma_k(\nu) = R_{k+1}(\mathfrak{m}_\nu) + 
P_k\bigl( R_2(\mathfrak{m}_\nu), \cdots, R_{k-1}(\mathfrak{m}_\nu)\bigr). 
\end{equation}
Here $P_k$ is a polynomial in $(k-2)$ variables. 
The part $P_k(\cdots)$ in \eqref{eq:2-3-5} has weight degree $\leqq k-1$. 
Note that $\Sigma_2(\nu) = R_3(\mathfrak{m}_\nu)$. 
See e.g. \cite{Bia03} for the Kerov polynomial. 
Hence we have for $k\in\mathbb{N}$ 
\begin{equation}\label{eq:2-3-6}
\Sigma_k(\nu) = R_{k+1}(\mathfrak{m}_\nu) + 
\sum_{\tau\in\mathbb{Y}: \, |\tau|+l(\tau)\leqq k-1} c_\tau \Sigma_\tau (\nu), 
\qquad c_\tau\in\mathbb{R}.
\end{equation}

\section{Evolution of averaged sizes}

Consider our continuous time process $(X_s^{(n)})_{s\geqq 0}$ defined by \eqref{eq:1-14}. 
Keeping \eqref{eq:1-15} and \eqref{eq:2-2-3} in mind, set 
\begin{equation}\label{eq:3-1}
a(k,n,s) = \sum_{j=0}^\infty \bigl(1- \frac{k}{n}\bigr)^j \int_{[0, s]} 
\psi\bigl( (s-u, \infty)\bigr) \psi^{\ast j}(du), \qquad k\in\mathbb{N}.
\end{equation}
Note $a(k,n,0) =1$ by $\psi(\{0\}) =0$. 
\eqref{eq:3-1} remains valid for $k=0$ if we set $a(0,n,s) =1$. 

\begin{proposition}\label{prop:3-1} 
The averaged size of each entry $X_s^{(n)\zeta}$ satisfies 
\[ 
\mathbb{E} \Bigl[ \frac{|X_s^{(n)\zeta}|}{n} \Bigr] = 
\bigl( 1- a(1,n,s)\bigr) \frac{(\dim\zeta)^2}{|T|} + 
a(1,n,s) \mathbb{E} \Bigl[ \frac{|X_0^{(n)\zeta}|}{n} \Bigr], \qquad \zeta\in\widehat{T}.
\] 
\end{proposition}
\textit{Proof} \ Let $M_s^{(n)}$ be the distribution of $X_s^{(n)}$, regarded as a row vector. 
\eqref{eq:1-15} yields 
\begin{align}
M_s^{(n)}(\mu) &= \sum_\lambda \mathbb{P} (X_s^{(n)} =\mu \,|\, 
X_0^{(n)} =\lambda) \,\mathbb{P}(X_0^{(n)} =\lambda) \notag \\ 
&= \sum_{j=0}^\infty (M_0^{(n)} P^{(n)j})_\mu \int_{[0, s]} 
\psi\bigl( (s-u, \infty)\bigr) \psi^{\ast j}(du), \qquad \mu\in\mathbb{Y}_n(\widehat{T}). 
\label{eq:3-3}
\end{align}
\eqref{eq:2-1-5} yields 
\begin{equation}\label{eq:3-4}
\frac{\chi^\mu_{\iota_{n,1}(1)_\theta}}{\dim\mu} = \sum_{\zeta\in\widehat{T}} 
\frac{\chi^\zeta_\theta}{\dim\zeta} \, \frac{|\mu^\zeta|}{n}.
\end{equation}
Applying \eqref{eq:2-2-3} to $(1, \theta)$-cycle ($\theta\neq\{e_T\}$), we have 
\begin{equation}\label{eq:3-5}
P^{(n) j} \Bigl( \frac{\chi^\mu_{\iota_{n,1}(1)_\theta}}{\dim\mu}
\Bigr)_{\mu\in\mathbb{Y}_n(\widehat{T})} = \bigl( 1- \frac{1}{n}\bigr)^j 
\Bigl( \frac{\chi^\mu_{\iota_{n,1}(1)_\theta}}{\dim\mu}
\Bigr)_{\mu\in\mathbb{Y}_n(\widehat{T})}, \qquad \theta\neq\{e_T\}.
\end{equation}
Combining \eqref{eq:3-3} with \eqref{eq:3-5}, we have 
\begin{equation}\label{eq:3-6}
\sum_\mu M_s^{(n)}(\mu) \frac{\chi^\mu_{\iota_{n,1}(1)_\theta}}{\dim\mu} 
= a(1,n,s) \sum_\mu M_0^{(n)}(\mu) \frac{\chi^\mu_{\iota_{n,1}(1)_\theta}}{\dim\mu} 
\qquad \theta\neq\{e_T\}.
\end{equation}
Putting \eqref{eq:3-4} into \eqref{eq:3-6} and adding the first column, we have  in matrix forms 
\begin{equation}\label{eq:3-7}
\Bigl( \mathbb{E} \Bigl[ \frac{|X_s^{(n)\zeta}|}{n}\Bigr] \Bigr)_\zeta 
\Bigl( \frac{\chi^\zeta_\theta}{\dim\zeta} \Bigr)_{\zeta, \theta} = 
\Bigl( \mathbb{E} \Bigl[ \frac{|X_0^{(n)\zeta}|}{n}\Bigr] \Bigr)_\zeta 
\begin{pmatrix} 1&\ \\ \vdots & a(1,n,s) \dfrac{\chi^\zeta_\theta}{\dim\zeta} \\ 
1 &\ \end{pmatrix}_{\zeta, \,\theta\neq\{e_T\}}
\end{equation}
Since 
\[ 
\Bigl( 1 \ O\Bigr) \Bigl( \frac{\chi^\zeta_\theta}{\dim\zeta} \Bigr)^{-1} 
= \Bigl( \frac{(\dim\zeta)^2}{|T|} \Bigr)_{\theta, \zeta}
\] 
holds, \eqref{eq:3-7} yields as row vectors 
\[ 
\Bigl( \mathbb{E} \Bigl[ \frac{|X_s^{(n)\zeta}|}{n}\Bigr] \Bigr)_\zeta 
= \bigl( (1-a(1,n,s)) \frac{(\dim\zeta)^2}{|T|}\bigr)_\zeta + 
a(1,n,s) \Bigl( \mathbb{E} \Bigl[ \frac{|X_0^{(n)\zeta}|}{n}\Bigr] \Bigr)_\zeta. 
\] 
as desired. 
\hfill $\blacksquare$

\bigskip

Let us turn to the rescale of $s = t \tau_n$ as \eqref{eq:1-17}. 
We recall several cases where the limit of $a(k,n, t\tau_n)$ in \eqref{eq:3-1} exists 
as $n\to\infty$. 

\begin{proposition}\label{prop:3-2} 
(1) Let $\psi$ have mean $m$ and its characteristic function $\varphi$ satisfy 
\[ 
\int_{|u|\geqq \delta} \Bigl| \frac{\varphi(u)}{u}\Bigr| du < \infty.
\] 
Set $\tau_n =n$. 
Then, we have 
\begin{equation}\label{eq:3-11}
\lim_{n\to\infty} a(k,n,tn) = e^{-kt/m}, \qquad k\in\mathbb{N}, \ t\geqq 0.
\end{equation}
(2) Let $\psi$ be the one-sided stable distribution of exponent $\alpha\in (0,1)$ i.e. 
its characteristic function $\varphi$ be given by 
\[ 
\varphi(u) = e^{-|u|^\alpha ( 1- i\tan(\pi \alpha/2) \mathrm{sgn}\,u)}, \qquad u\in\mathbb{R}.
\] 
Set $\tau_n = n^{1/\alpha}$. 
Then, we have 
\begin{equation}\label{eq:3-13}
\lim_{n\to\infty} a(k,n,tn^{1/\alpha}) = \frac{\sin \pi\alpha}{\pi\alpha} \int_0^\infty 
\frac{e^{-t(ku \cos (\pi\alpha/2))^{1/\alpha}}}{u^2+2u\cos\pi\alpha +1}\, du, 
\qquad k\in\mathbb{N}, \ t\geqq 0.
\end{equation}
\end{proposition}
\textit{Proof} \ 
See Propositions 2.1 and 2.2 in \cite{Hor20}
\footnote{We set $k\geqq 2$ in those propositions in \cite{Hor20} since the case of $k=1$ 
was unnecessary. 
However, those propositions themselves are valid for $k=1$.}. 
\hfill 
$\blacksquare$ 

\medskip

If pausing time distribution $\psi$ satisfies (1) or (2) of Proposition\ref{prop:3-2}, set 
\begin{equation}\label{eq:3-15}
a_k(t) = \lim_{n\to\infty} a(k, n, t\tau_n), \qquad k\in\mathbb{N}, \ t\geqq 0
\end{equation}
for \eqref{eq:3-11} or \eqref{eq:3-13}. 

\begin{theorem}\label{th:3-3}
For random multi-diagrams $(X_{t\tau_n}^{(n)})$, assume that the initial ensemble satisfies 
\begin{equation}\label{eq:3-14}
\lim_{n\to\infty} \mathbb{E} \Bigl[ \frac{|X_0^{(n)\zeta}|}{n}\Bigr] 
= n^\zeta, \qquad \zeta\in\widehat{T}
\end{equation}
and that pausing time distribution $\psi$ satisfies (1) or (2) of Proposition\ref{prop:3-2}. 
Then, the averaged size of each entry at time $t\geqq 0$ is given by 
\begin{equation}\label{eq:3-16}
\lim_{n\to\infty} \mathbb{E} \Bigl[ \frac{|X_{t\tau_n}^{(n)\zeta}|}{n}\Bigr] 
= (1-a_1(t)) \frac{(\dim\zeta)^2}{|T|} + a_1(t) n^\zeta, \qquad \zeta\in\widehat{T}.
\end{equation}
In particular, the averaged sizes of the entries in the limit of $t\to\infty$ are distributed 
according to the Plancherel measure of $T$.
\end{theorem}
\textit{Proof} \ 
\eqref{eq:3-16} follows from Propositions \ref{prop:3-1} and \ref{prop:3-2}. 
Since \eqref{eq:3-11} and \eqref{eq:3-13} yield 
\[ 
\lim_{t\to\infty} a_1(t) =0
\] 
in \eqref{eq:3-15}, the final assertion is valid. 
\hfill $\blacksquare$

\bigskip

\noindent\textit{Example} 
We mention an obvious example for Theorem\ref{th:3-3}. 
For arbitrary $(n^\zeta)_{\zeta\in\widehat{T}}$ satisfying $n^\zeta >0$ and 
$\sum_\zeta  n^\zeta =1$, 
let a multiple continuous diagram $(\omega^\zeta)_{\zeta\in\widehat{T}}$ be given 
so that 
\[ 
\omega^\zeta\in \mathbb{D}, \qquad 
\int_\mathbb{R} (\omega^\zeta(x)-|x|)\, dx = 2 n^\zeta.
\] 
A continuous diagram is approximated uniformly by rescaled Young diagrams \eqref{eq:1-16}. 
Namely, we can take 
$\lambda_n = (\lambda_n^\zeta)_{\zeta\in\widehat{T}} \in \mathbb{Y}_n(\widehat{T})$ 
such that 
\[ 
\lim_{n\to\infty} (\lambda_n^\zeta)^{\sqrt{n}}(x) = \omega^\zeta(x) \quad 
\text{uniformly in } x, \qquad \zeta\in\widehat{T}.
\] 
Initial ensemble $M_0^{(n)} = \delta_{\lambda_n}$ satisfies \eqref{eq:3-14}. 
Indeed, for any $\zeta\in\widehat{T}$, 
\[ 
\mathbb{E} \Bigl[ \frac{|X_0^{(n)\zeta}|}{n}\Bigr] = \frac{|\lambda_n^\zeta|}{n} 
= \frac{1}{2} \int_\mathbb{R} \bigl( (\lambda_n^\zeta)^{\sqrt{n}}(x)-|x|\bigr)\, dx \ 
\xrightarrow[n\to\infty] \ \frac{1}{2}\int_\mathbb{R} (\omega^\zeta(x)-|x|)\, dx = n^\zeta.
\] 

\section{Evolution of averaged limit shapes}

Theorem\ref{th:3-3} deals only with the size (number of boxes) of a Young diagram, 
which is the variance (the second cumulant) of its transition measure. 
In order to characterize the (averaged) limit shape, we need to know all the (free) cumulants 
of higher degrees. 

A probability $M$ on $\mathbb{Y}_n(\widehat{T})$ determines a function $f$ on 
$\mathfrak{S}_n(T)$ by 
\begin{equation}\label{eq:4-1}
f(x) = \sum_{\lambda\in\mathbb{Y}_n(\widehat{T})} M(\lambda) 
\frac{\chi^\lambda(x)}{\dim\lambda}, \qquad x\in \mathfrak{S}_n(T).
\end{equation}
\eqref{eq:4-1} gives a bijective correspondence between 
$\mathcal{M}(\mathbb{Y}_n(\widehat{T}))$, the probabilities on $\mathbb{Y}_n(\widehat{T})$, 
and $\mathcal{K}(\mathfrak{S}_n(T))$, the positive-definite central normalized functions 
on $\mathfrak{S}_n(T)$. 
When we consider assumptions on an initial ensemble, it is often more convenient to state them 
in terms of $\mathcal{K}(\mathfrak{S}_n(T))$ than 
$\mathcal{M}(\mathbb{Y}_n(\widehat{T}))$. 
In fact, we later discuss the ensembles obtained by restricting an extremal element of 
$\mathcal{K}(\mathfrak{S}_\infty (T))$ onto $\mathfrak{S}_n(T)$. 

In our process $(X_s^{(n)})_{s\geqq 0}$, let 
$M_s^{(n)}\in \mathcal{M}(\mathbb{Y}_n(\widehat{T}))$ 
be the distribution at $s$, namely 
$M_s^{(n)}(\lambda) = \mathbb{P}( X_s^{(n)} =\lambda)$, and 
$f_s^{(n)}\in \mathcal{K}(\mathfrak{S}_n(T))$ the corresponding positive-definite 
function of \eqref{eq:4-1}. 
Recalling \eqref{eq:3-4}, we see from \eqref{eq:4-1} 
the assumption for initial ensembles \eqref{eq:3-14} in Theorem\ref{th:3-3} 
\[ 
n^\zeta = \lim_{n\to\infty} \sum_{\mu\in\mathbb{Y}_n(\widehat{T})} 
M_0^{(n)}(\mu) 
\frac{|\mu^\zeta|}{n}, \qquad \zeta\in \widehat{T}
\] 
is equivalent to 
\begin{equation}\label{eq:4-4}
\lim_{n\to\infty} f_0^{(n)} \bigl( \iota_{n,1}(1)_\theta\bigr) = 
\sum_{\zeta\in\widehat{T}} \frac{\chi^\zeta_\theta}{\dim\zeta}\, n^\zeta, 
\qquad \theta\in [T].
\end{equation}
It is then natural to extend \eqref{eq:4-4} to sharpen Theorem\ref{th:3-3} 
(see \eqref{eq:4-5} below). 

For $\rho = (\rho_\theta)_{\theta\in [T]}$, set 
$l(\rho) = \sum_{\theta\in [T]} l(\rho_\theta)$. 
The $j$th free cumulant of probability $\mathfrak{m}$ on $\mathbb{R}$ is denoted by 
$R_j(\mathfrak{m})$ ($j\in\mathbb{N}$). 

\begin{theorem}\label{th:4-1}
For random multi-diagrams $(X_{t\tau_n}^{(n)})$, assume that the initial ensemble satisfies 
\begin{align}
\lim_{n\to\infty} n^{\frac{k-1}{2}} f_0^{(n)} \bigl( \iota_{n,k}(k)_\theta\bigr) 
&= \gamma_{k+1}^\theta, \qquad k\in\mathbb{N}, \ \theta\in [T], 
\label{eq:4-5} \\ 
n^{\frac{k-l(\rho)}{2}} f_0^{(n)}( \iota_{n,k} \,\rho) &= O(1) \quad (n\to\infty), 
\qquad k\in\mathbb{N}, \ \rho\in\mathbb{Y}_k([T]). 
\label{eq:4-6}
\end{align}
Let the pausing time distribution $\psi$ satisfy (1) or (2) of Proposition\ref{prop:3-2} and 
set $a_k(t)$ by \eqref{eq:3-15} for \eqref{eq:3-11} or \eqref{eq:3-13}. 
Then, we have for $t\geqq 0$ 
\begin{align}
&\lim_{n\to\infty} n^{\frac{k-1}{2}} f_{t\tau_n}^{(n)} \bigl( \iota_{n,k}(k)_\theta\bigr) 
= a_k(t)\, \gamma_{k+1}^\theta, \qquad 
k\geqq 2, \ \theta\in [T], \text{ or }  k=1, \ \theta\neq \{e_T\}, 
\notag \\ 
&\lim_{n\to\infty} \mathbb{E}_{M_{t\tau_n}^{(n)}}\bigl[ 
R_{k+1}\bigl( \mathfrak{m}_{(\lambda^\zeta)^{\sqrt{n}}}\bigr) \bigr] 
= a_k(t) \sum_{\theta\in [T]} \gamma_{k+1}^\theta 
\frac{|C_\theta| \overline{\chi}^\zeta_\theta (\dim\zeta)^k}{|T|}, \qquad 
k\geqq 2, \ \zeta\in \widehat{T}. 
\label{eq:4-8}
\end{align}
\end{theorem}

\noindent\textit{Remark} \ 
We mention several comments on Theorem\ref{th:4-1}. 

\noindent (i) 
For each $k\in\mathbb{N}$, $(\gamma_{k+1}^\theta)_{\theta\in [T]}$ in \eqref{eq:4-5} 
detemines $(R_{k+1}^\zeta)_{\zeta\in\widehat{T}}$ (and vice versa) by 
\begin{equation}\label{eq:4-9}
\gamma_{k+1}^\theta = \sum_{\zeta\in\widehat{T}} 
\frac{\chi^\zeta_\theta}{(\dim\zeta)^k} R_{k+1}^\zeta \  (\theta\in [T]), \quad 
R_{k+1}^\zeta = \sum_{\theta\in [T]} \gamma_{k+1}^\theta \frac{|C_\theta|}{|T|}
\overline{\chi}^\zeta_\theta (\dim\zeta)^k \  (\zeta\in\widehat{T}), 
\end{equation}
where $C_\theta$ is the conjugacy class of $T$ labeled by $\theta$. 
In particular, we have $R_2^\zeta = n^\zeta$ by comparing \eqref{eq:4-9} 
with \eqref{eq:4-4}.

\noindent (ii) 
\eqref{eq:4-8} does not hold if $k=1$. 
For $k=1$, the LHS of \eqref{eq:4-8} is given by \eqref{eq:3-16}. 
On the other hand, the RHS is just $a_1(t) n^\zeta$ by \eqref{eq:4-9}. 

\noindent (iii) 
\eqref{eq:4-9} tells us that \eqref{eq:4-8} is equivalent to 
\begin{equation}\label{eq:4-10}
\lim_{n\to\infty} \mathbb{E}_{M_{t\tau_n}^{(n)}}\bigl[ 
R_{k+1}\bigl( \mathfrak{m}_{(\lambda^\zeta)^{\sqrt{n}}}\bigr) \bigr] 
= a_k(t) R_{k+1}^\zeta.
\end{equation}
Combining \eqref{eq:3-16}, 
we thus obtain a sequence $(R_j^\zeta(t))$ for each $\zeta\in\widehat{T}$ determined by 
\begin{equation}\label{eq:4-11}
R_1^\zeta(t) =0, \quad 
R_2^\zeta(t) = (1-a_1(t)) \frac{(\dim\zeta)^2}{|T|} + a_1(t) R_2^\zeta, \quad 
R_{k+1}^\zeta(t) = a_k(t) R_{k+1}^\zeta \ (k\geqq 2).
\end{equation}
If $(R_j^\zeta(t))_{j\in\mathbb{N}}$ is the sequence of free cumulants of a probability 
$\mathfrak{m}^\zeta(t)$ on $\mathbb{R}$, we can regard the multiple continuous diagram 
$(\omega^\zeta(t))_{\zeta\in\widehat{T}}$ determined by the 
transition measures $(\mathfrak{m}^\zeta(t))_{\zeta\in\widehat{T}}$ as the 
averaged limit shape at time $t$ for our random multi-diagrams. 
Later we discuss an example realizing this situation. 

\noindent (iv) 
After the limit of $t\to\infty$, \eqref{eq:4-11} turns into 
\begin{equation}\label{eq:4-12}
R_j^\zeta (\infty) =0 \quad (j\neq 2), \qquad 
R_2^\zeta (\infty) = \frac{(\dim\zeta)^2}{|T|}, 
\end{equation}
as easily seen from \eqref{eq:3-11} and \eqref{eq:3-13}. 
\eqref{eq:4-12} is the free cumulant sequence of the semi-circle distribution with 
mean $0$ and variance $(\dim\zeta)^2 /|T|$, which is the transition measure of the 
Vershik--Kerov--Logan--Shepp curve rescaled by $\dim\zeta /\sqrt{|T|}$. 

\noindent (v) 
Concerning the assumption \eqref{eq:4-6}, let 
$\rho = (\rho_\theta)\in \mathbb{Y}_k([T])$ with $m_1(\rho_{\{e_T\}}) =0$ and 
$\rho^\prime = (\rho_{\{e_T\}}\sqcup (1^j), \rho_{\theta})_{\theta\neq\{e_T\}}
\in \mathbb{Y}_{k+j}([T])$. 
Then, $l(\rho^\prime) = l(\rho)+j$ and 
$\iota_{n,k}\,\rho = \iota_{n,k+j}\, \rho^\prime$ hold. 
Hence \eqref{eq:4-6} for $\rho^\prime$ follows from \eqref{eq:4-6} for $\rho$. 
In other words, \eqref{eq:4-6} can be restricted to the case of $m_1(\rho_{\{e_T\}}) =0$ 
without loss of generality. 

\bigskip

\noindent\textit{Proof of Theorem\ref{th:4-1}} \ 
First we note 
\[ 
f_s^{(n)}( \iota_{n,k}\rho) = a(k-m_1(\rho_{\{e_T\}}), n, s) f_0^{(n)}(\iota_{n,k}\rho) 
\] 
holds for $k\in\mathbb{N}$ and $\rho \in\mathbb{Y}_k([T])$. 
In fact, \eqref{eq:3-3} and \eqref{eq:2-2-2} yield 
\begin{align*}
&f_s^{(n)}( \iota_{n,k}\rho) = 
\sum_\mu M_s^{(n)}(\mu) \frac{\chi^\mu_{\iota_{n,k}\rho}}{\dim\mu} \\ 
&= \sum_\mu \sum_{j=0}^\infty (M_0^{(n)} P^{(n) j})_\mu 
\frac{\chi^\mu_{\iota_{n,k}\rho}}{\dim\mu} 
\int_{[0, s]} \psi\bigl( (s-u, \infty)\bigr) \psi^{\ast j} (du)  \\ 
&= \sum_{j=0}^\infty \Bigl( 1- \frac{k-m_1(\rho_{\{e_T\}})}{n}\Bigr)^j 
\Bigl( \int_{[0, s]} \psi\bigl( (s-u, \infty)\bigr) \psi^{\ast j} (du) \Bigr)
\mathbb{E}_{M_0^{(n)}}\Bigl[ \frac{\chi^\mu_{\iota_{n,k}\rho}}{\dim\mu} \Bigr]. 
\end{align*}
Then, for $t\geqq 0$, we have 
\begin{equation}\label{eq:4-15}
n^{\frac{k-1}{2}} f_{t\tau_n}^{(n)} (\iota_{n,k} (k)_\theta) 
= a(k,n, t\tau_n) n^{\frac{k-1}{2}} f_0^{(n)} (\iota_{n,k}\rho) 
\ \xrightarrow[n\to\infty] \ a_k(t) \gamma_{k+1}^\theta
\end{equation}
for $k\geqq 2$ or $k=1$ and $\theta\neq\{e_T\}$ by \eqref{eq:4-5} and \eqref{eq:3-15}, and 
\begin{equation}\label{eq:4-16}
n^{\frac{k-l(\rho)}{2}} f_{t\tau_n}^{(n)}(\iota_{n,k}\rho) = O(1) \quad (n\to\infty), 
\qquad \rho\in\mathbb{Y}_k([T])
\end{equation}
by \eqref{eq:4-6} and \eqref{eq:3-15}. 
Since the LHS of \eqref{eq:4-15} is 
\[ 
\frac{n^k}{n^{\downarrow k}} \sum_{\zeta\in\widehat{T}} 
\frac{\chi^\zeta_\theta}{(\dim\zeta)^k} \mathbb{E}_{M_{t\tau_n}^{(n)}} 
\Bigl[ n^{-\frac{k+1}{2}} \Sigma_k(\lambda^\zeta) \Bigr] 
\] 
by \eqref{eq:2-1-5}, we see from \eqref{eq:4-15} and \eqref{eq:4-9} the convergence 
\begin{equation}\label{eq:4-18}
\lim_{n\to\infty} \mathbb{E}_{M_{t\tau_n}^{(n)}} 
\Bigl[ n^{-\frac{k+1}{2}} \Sigma_k(\lambda^\zeta) \Bigr] 
= a_k(t) R_{k+1}^\zeta.
\end{equation}
\eqref{eq:2-3-6} enables us to replace $\Sigma_k(\lambda^\zeta)$ by 
$R_{k+1}(\mathfrak{m}_{\lambda^\zeta})$ in \eqref{eq:4-18}. 
To be precise, that follows from Lemma\ref{lem:4-2}. 
Then, noting 
\[ 
n^{-\frac{k+1}{2}} R_{k+1}(\mathfrak{m}_{\lambda^\zeta}) = 
R_{k+1}(\mathfrak{m}_{(\lambda^\zeta)^{\sqrt{n}}}), 
\] 
we finally have \eqref{eq:4-10}. 
\hfill $\blacksquare$

\begin{lemma}\label{lem:4-2}
Condition \eqref{eq:4-16} for $k\in\mathbb{N}$ yields 
\[ 
\lim_{n\to\infty} \mathbb{E}_{M_{t\tau_n}^{(n)}} 
\Bigl[ n^{-\frac{k+1}{2}} \Sigma_\tau (\lambda^\zeta) \Bigr] =0
\] 
for $\zeta\in\widehat{T}$ and $\tau\in\mathbb{Y}$ such that 
$|\tau|+l(\tau) \leqq k-1$.
\end{lemma}
\textit{Proof} \ 
Let $\tau\in\mathbb{Y}$ be so fixed and $\mathrm{rows}(\tau)$ denote the set 
of rows of $\tau$. 
For $s : \mathrm{rows}(\tau) \longrightarrow [T]$, each $s^{-1}(\theta)$ is 
regarded as a Young diagram. 
Obviously
\footnote{Since all rows are distinguished, $(s^{-1}(\theta)) = (s^{\prime -1}(\theta))$ 
can occur even if $s\neq s^\prime$.}, 
\[ 
\bigl( s^{-1}(\theta) \bigr)_{\theta\in [T]} \in \mathbb{Y}_{|\tau|}([T]), \quad 
l(\tau) = l((s^{-1}(\theta))) \quad \text{and} \quad 
\mathrm{rows}\bigl( (s^{-1}(\theta)) \bigr) = \mathrm{rows}(\tau).
\] 
Since \eqref{eq:2-1-4} gives 
\[ 
\frac{\chi^\lambda_{\iota_{n,|\tau|}(s^{-1}(\theta))_{\theta\in [T]}}}{\dim\lambda} 
= \frac{1}{n^{\downarrow |\tau|}} \!\!\sum_{r:\, \mathrm{rows}(\tau) \to \widehat{T}} 
\prod_{\zeta\in\widehat{T}} \Bigl\{ 
\frac{\Sigma_{r^{-1}(\zeta)}(\lambda^\zeta)}{(\dim\zeta)^{|r^{-1}(\zeta)|}} 
\prod_{\theta\in [T]} (\chi^\zeta_\theta)^{l(r^{-1}(\zeta)\cap s^{-1}(\theta))}\Bigr\}, 
\] 
we have 
\begin{multline}\label{eq:4-23}
n^{\frac{|\tau|-l(\tau)}{2}} f_{t\tau_n}^{(n)} 
\bigl(\iota_{n,|\tau|}(s^{-1}(\theta))_{\theta\in [T]}\bigr) \\ 
= \frac{n^{|\tau|}}{n^{\downarrow |\tau|}} \!\!\sum_{r:\, \mathrm{rows}(\tau) \to \widehat{T}} 
\Bigl( \prod_{\zeta\in\widehat{T}} \prod_{\theta\in [T]} 
\frac{(\chi^\zeta_\theta)^{l(r^{-1}(\zeta)\cap s^{-1}(\theta))}}
{(\dim\zeta)^{|r^{-1}(\zeta)|}} \Bigr) 
\frac{1}{n^{\frac{|\tau|+l(\tau)}{2}}} \mathbb{E}_{M_{t\tau_n}^{(n)}} \Bigl[ 
\prod_{\zeta\in\widehat{T}}\Sigma_{r^{-1}(\zeta)}(\lambda^\zeta) \Bigr]. 
\end{multline}
In \eqref{eq:4-23}, the $\widehat{T}^{l(\tau)} \times [T]^{l(\tau)}$ matrix 
\[ 
A = \Bigl( \prod_{\zeta\in\widehat{T}} \prod_{\theta\in [T]} 
(\chi^\zeta_\theta)^{l(r^{-1}(\zeta)\cap s^{-1}(\theta))}\Bigr)_{r,s} 
\] 
is nothing but the $l(\tau)$-fold tensor product of $(\chi^\zeta_\theta)$ and hence 
invertible. 
Then \eqref{eq:4-23} yields (as row vectors) 
\begin{multline}\label{eq:4-25}
\Bigl( n^{\frac{|\tau|-l(\tau)}{2}} f_{t\tau_n}^{(n)} 
\bigl(\iota_{n,|\tau|}(s^{-1}(\theta))_{\theta\in [T]}\bigr) \Bigr)_s A^{-1} 
\mathrm{diag}\Bigl( \prod_{\zeta\in\widehat{T}} (\dim\zeta)^{|r^{-1}(\zeta)|}\Bigr)_r \\ 
= \frac{n^{|\tau|}}{n^{\downarrow |\tau|}} \Bigl( 
\frac{1}{n^{\frac{|\tau|+l(\tau)}{2}}} \mathbb{E}_{M_{t\tau_n}^{(n)}} \Bigl[ 
\prod_{\zeta\in\widehat{T}}\Sigma_{r^{-1}(\zeta)}(\lambda^\zeta) \Bigr] \Bigr)_r.
\end{multline}
In the RHS of \eqref{eq:4-25}, we pick up the $r$ entry such that $r^{-1}(\zeta) = \tau$ 
(i.e. $r \equiv \zeta$) holds for initially given $\zeta\in\widehat{T}$. 
\eqref{eq:4-16} implies the LHS of \eqref{eq:4-25} is $O(1)$. 
Hence we have 
\[ 
 \mathbb{E}_{M_{t\tau_n}^{(n)}} \bigl[ \Sigma_\tau(\lambda^\zeta) \bigr] 
= O\bigl(n^{\frac{|\tau|+l(\tau)}{2}}\bigr). 
\] 
This completes the proof.
\hfill $\blacksquare$

\medskip

We got \eqref{eq:3-11} and \eqref{eq:3-13} as expressions of $a_k(t)$ according to 
the conditions for pausing time distributions. 
Let us make further consideration in these cases separately. 
We refer to \cite{VoDyNi92} and \cite{NiSp06} for notions in free probability theory. 

\begin{corollary}\label{cor:4-3}
Consider $a_k(t) = e^{-kt/m}$ ($k\in\mathbb{N}$) as \eqref{eq:3-11}, and assume 
$(R^\zeta_j)_{j\in\mathbb{N}}$ determined by \eqref{eq:4-9} 
(with $R^\zeta_1 =0$ and $R^\zeta_2 = n^\zeta$) 
to be the free cumulant sequence of a probability $\mathfrak{m}^\zeta$ on $\mathbb{R}$. 
Then, $(R^\zeta_j(t))_{j\in\mathbb{N}}$ determined by \eqref{eq:4-11} is realized as 
the free cumulant sequence of the probability 
\begin{equation}\label{eq:4-27}
\mathfrak{m}^\zeta (t) = (\mathfrak{m}^\zeta)_{e^{-t/m}} \boxplus 
(\gamma^\zeta)_{1-e^{-t/m}}.
\end{equation}
Here $\gamma^\zeta$ is the semi-circle distribution with mean $0$ and variance 
$(\dim\zeta)^2/|T|$, $\boxplus$ denotes the free (additive) convolution, and 
subscription $_c$ denotes the free compression by rank $c$ projection. 
Voiculescu's $R$-transform of $\mathfrak{m}^\zeta (t)$ is given by 
\begin{equation}\label{eq:4-28}
R^\zeta (t, w) = \sum_{k=0}^\infty R_{k+1}\bigl( \mathfrak{m}^\zeta (t)\bigr) w^k 
= (1-e^{-\frac{t}{m}}) \frac{(\dim\zeta)^2}{|T|} w + 
R^\zeta (0, e^{-\frac{t}{m}} w).
\end{equation}
Furthermore, the Stieltjes transform $G^\zeta(t, z)$ of $\mathfrak{m}^\zeta (t)$ 
satisfies PDE 
\begin{equation}\label{eq:4-29}
\frac{\partial G^\zeta}{\partial t} = \Bigl( \frac{1}{m G^\zeta} 
-\frac{(\dim\zeta)^2}{|T|}\, \frac{G^\zeta}{m} \Bigr) \frac{\partial G^\zeta}{\partial z} 
+ \frac{G^\zeta}{m}.
\end{equation}
\end{corollary}
\textit{Proof} proceeds along quite similar computations to our previous results in 
\cite{Hor15} and \cite{Hor20}, and is omitted. 
\eqref{eq:4-27} rephrases \eqref{eq:4-11}. 
\hfill $\blacksquare$

\medskip

The case of \eqref{eq:3-13} is less easy. 
We discuss only the case of $\alpha =1/2$ here, so let 
\begin{equation}\label{eq:4-30}
a_k(t) = \frac{2}{\pi} \int_0^\infty \frac{e^{-tk^2 u^2/2}}{u^2+1} \, du, \qquad 
k\in\mathbb{N}, \quad t\geqq 0.
\end{equation}

\begin{proposition}\label{prop:4-4}
Let $a_k(t)$ be as \eqref{eq:4-30}, and assume 
$(R^\zeta_j)_{j\in\mathbb{N}}$ determined by \eqref{eq:4-9} 
(with $R^\zeta_1 =0$ and $R^\zeta_2 = n^\zeta$) 
to be the free cumulant sequence of a compactly supported $\boxplus$-infinitely divisible 
probability $\mathfrak{m}^\zeta$ on $\mathbb{R}$ 
\footnote{Later we mention \eqref{eq:4-48}, for example.}. 
Then, $(R^\zeta_j(t))_{j\in\mathbb{N}}$ determined by \eqref{eq:4-11} is also 
the free cumulant sequence of compactly supported $\boxplus$-infinitely divisible probability 
$\mathfrak{m}^\zeta(t)$ on $\mathbb{R}$. 
Voiculescu's $R$-transform of $\mathfrak{m}^\zeta (t)$ is given by 
\begin{equation}\label{eq:4-31}
R^\zeta (t, w) = \int_\mathbb{R} \frac{1}{\sqrt{2\pi}} e^{-u^2/2}\Bigl\{ 
(1-e^{-\sqrt{t}|u|}) \frac{(\dim\zeta)^2}{|T|} w + 
R^\zeta (0, e^{-\sqrt{t}|u|}w) \Bigr\} du.
\end{equation}
\end{proposition}
\textit{Proof} \ 
We show $(R^\zeta_j(t))_{j\in\mathbb{N}}$ is conditionally positive-definite. 
See \cite{NiSp06} for the relation between this notion and $\boxplus$-infinite divisibility. 
In fact, 
for $r\in\mathbb{N}$ and $\alpha_j \in\mathbb{C}$, we have 
\begin{equation}\label{eq:4-32}
\sum_{j,k=1}^r \alpha_j \overline{\alpha_k} R^\zeta_{j+k} (t) 
= |\alpha_1|^2 (1-a_1(t)) \frac{(\dim\zeta)^2}{|T|} + 
\sum_{j,k=1}^r \alpha_j \overline{\alpha_k} a_{j+k-1}(t) R^\zeta_{j+k}. 
\end{equation}
Since \eqref{eq:4-30} is rewritten as 
\begin{equation}\label{eq:4-33}
a_k(t) = \frac{1}{\sqrt{2\pi t}} \int_\mathbb{R} e^{-x^2/(2t)} e^{-k|x|} dx
\end{equation}
through Fourier transform, the second sum of \eqref{eq:4-32} is equal to 
\[ 
\frac{1}{\sqrt{2\pi t}} \int_\mathbb{R} e^{-\frac{x^2}{2t}+|x|} 
\Bigl( \sum_{j,k=1}^r (\alpha_j e^{-j|x|}) ( \overline{\alpha_k e^{-k|x|}}) 
R^\zeta_{j+k} \Bigr) dx \geqq 0.
\] 
Hence we have \eqref{eq:4-32} $\geqq 0$. 
Corresponding probability $\mathfrak{m}^\zeta(t)$ is constructed by following 
the prescription (13.10) in \cite{NiSp06}. 
The $R$-transform $R^\zeta (t, w)$ of $\mathfrak{m}^\zeta(t)$ is computed from 
\eqref{eq:4-11} and \eqref{eq:4-33} to be \eqref{eq:4-31}.
\hfill $\blacksquare$

\medskip

\noindent\textit{Remark} \ 
(i) We note that $\mathfrak{m}^\zeta(t)$ obtained in Proposition\ref{prop:4-4} is 
the transition measure of an \textit{averaged} limit shape at time $t$. 
As we showed in \cite{Hor20}, we cannot expect concentration in the case of 
the pausing time distribution of Proposition\ref{prop:3-2} (2). 
In order to discuss concentration at a limit shape, the notion of approximate factorization 
property introduced by Biane \cite{Bia01} is useful. 
We develop such a concentrated limit shape for multi-diagrams in a separate note. 

\noindent (ii) Unfortunately, the above proof does not work for $\alpha\neq 1/2$. 

\noindent (iii) In Corollary\ref{cor:4-3} also, $\boxplus$-infinite divisibility is obviously 
inherited from $\mathfrak{m}^\zeta$ by $\mathfrak{m}^\zeta(t)$. 
The $R$-transform of a $\boxplus$-infinitely divisible probabiliry $\mathfrak{m}$ is 
expressed as 
\begin{equation}\label{eq:4-35}
R_\mathfrak{m}(w) = \int_\mathbb{R} \frac{w}{1-xw}\, \mathfrak{l} (dx) 
= G_{\mathfrak{l}}\bigl( \frac{1}{w}\bigr) 
\end{equation}
in terms of the Stieltjes transform $G_{\mathfrak{l}}$ of L\'evy measure $\mathfrak{l}$. 
Let $\mathfrak{l}_t^\zeta$ be the L\'evy measure of $\mathfrak{m}^\zeta (t)$ at $t\geqq 0$ in 
Collorary\ref{cor:4-3} or Proposition\ref{prop:4-4}. 
Then, \eqref{eq:4-28} yields 
\begin{equation}\label{eq:4-36}
\mathfrak{l}_t^\zeta (dx) = (1-e^{-t/m})\frac{(\dim\zeta)^2}{|T|}\, \delta_0(dx) + 
e^{-t/m} \,\mathfrak{l}_0^\zeta ( e^{t/m} dx),
\end{equation}
and \eqref{eq:4-31} yields 
\begin{equation}\label{eq:4-37}
\mathfrak{l}_t^\zeta (dx) = \int_\mathbb{R} \frac{e^{-u^2/2}}{\sqrt{2\pi}} 
(1-e^{-\sqrt{t}|u|})du \, \frac{(\dim\zeta)^2}{|T|}\, \delta_0 (dx) 
+ \Tilde{\mathfrak{l}}_t^\zeta (dx)
\end{equation}
respectively, where $\Tilde{\mathfrak{l}}_t^\zeta(dx)$ is the projection to the $x$ axis 
$\mathbb{R}$ of 
\[ 
\frac{1}{\sqrt{2\pi}} e^{-\frac{y^2}{2} - \sqrt{t}|y|} \,
\mathfrak{l}_0^\zeta( e^{\sqrt{t}|y|} dx) dy 
\] 
on the $xy$ plane $\mathbb{R}^2$.

\noindent (iv) 
Though we get an averaged limit shape from $\mathfrak{m}^\zeta$ through Markov 
transform in principle, concrete computation is often difficult. 
Partially, we have a PDE \eqref{eq:4-29} for the Stieltjes transform. 
However, we do not know functional equation governing the limit shape (like the ones 
obtained in \cite{FuSa10}) yet. 

\section{Example from characters of $\mathfrak{S}_\infty(T)$}

To conclude, let us mention an example of ensembles, a sequence of probabilities on 
$\mathbb{Y}_n(\widehat{T})$'s, realizing Theorem\ref{th:4-1}, Corollary\ref{cor:4-3} 
and Proposition\ref{prop:4-4}. 

The canonical action of the infinite symmetric group $\mathfrak{S}_\infty$ onto 
the restricted direct product of finite group, 
$T_\infty = \{ (t_j)_{j\in\mathbb{N}} \,|\, t_j\in T, \ t_j = e_T 
\text{ except finite $j$'s}\}$, defines the infinite wreath product 
$\mathfrak{S}_\infty(T)$. 
We refer to an extremal normalized central positive-definite function on 
$\mathfrak{S}_\infty(T)$ simply as a character. 
The characters of $\mathfrak{S}_\infty(T)$ are classified in \cite{Boy05} and 
\cite{HiHi05} when $T$ is finite. 
Here we use notations in \cite{HoHi14} though that paper is devoted mainly to the case of 
compact (continuous) $T$. 
The conjugacy classes of $\mathfrak{S}_\infty(T)$ are parametrized by 
\[ 
\mathbb{Y}([T]) = \bigl\{ \rho = (\rho_\theta)_{\theta\in [T]} \,\big|\, 
\rho_\theta\in\mathbb{Y}, \; \sum_{\theta\in [T]} |\rho_\theta| < \infty, \; 
m_1(\rho_{\{e_T\}})=0 \bigr\}.
\] 
The characters of $\mathfrak{S}_\infty(T)$ are parametrized by 
\begin{align}\label{eq:4-40}
\Delta = \bigl\{ (\alpha, \beta, c) \,\big|\, 
&\alpha = (\alpha_{\zeta, i})_{\zeta\in\widehat{T}, i\in\mathbb{N}}, \; 
\beta = (\beta_{\zeta, i})_{\zeta\in\widehat{T}, i\in\mathbb{N}}, \; 
c = (c_\zeta)_{\zeta\in\widehat{T}}, \ 
\alpha_{\zeta, 1}\geqq \alpha_{\zeta, 2}\geqq \cdots\geqq 0, \notag \\
&\beta_{\zeta, 1}\geqq \beta_{\zeta, 2}\geqq \cdots\geqq 0, \; c_\zeta \geqq 0, \ 
\sum_{i=1}^\infty (\alpha_{\zeta, i}+ \beta_{\zeta, i}) \leqq c_\zeta, \; 
\sum_{\zeta\in\widehat{T}} c_\zeta =1\bigr\}
\end{align}
with a character value 
\begin{equation}\label{eq:4-41}
f_\omega (\rho) = \prod_{j=1}^\infty \prod_{\theta\in [T]} \Bigl( 
\sum_{\zeta\in\widehat{T}} p^\zeta_j(\omega) \frac{\chi^\zeta_\theta}{(\dim\zeta)^j}
\Bigr)^{m_j(\rho_\theta)}, \qquad \omega\in\Delta, \ \rho\in\mathbb{Y}([T]), 
\end{equation}
where we set for $\omega = (\alpha, \beta, c)$ 
\begin{equation}\label{eq:4-42}
p^\zeta_j(\omega) = \begin{cases} \sum_{i=1}^\infty (\alpha_{\zeta, i}^{\ j} 
+ (-1)^{j-1} \beta_{\zeta, i}^{\ j}), & j\geqq 2 \\ c_\zeta, & j=1. \end{cases}
\end{equation}
A family of atomic measures $(\tau_\omega^\zeta)_{\zeta\in\widehat{T}}$ for 
$\omega= (\alpha, \beta, c) \in\Delta$ is defined by 
\begin{equation}\label{eq:4-43}
\tau_\omega^\zeta = \sum_{i=1}^\infty ( \alpha_{\zeta, i} \delta_{\alpha_{\zeta, i}}
+ \beta_{\zeta, i} \delta_{-\beta_{\zeta, i}}) + 
\bigl( c_\zeta - \sum_{i=1}^\infty (\alpha_{\zeta, i}+\beta_{\zeta, i})\bigr)\delta_0
\end{equation}
and called a Thoma measure. 
As its moment, \eqref{eq:4-42} and \eqref{eq:4-43} yield 
\begin{equation}\label{eq:4-44}
M_{k-1}(\tau_\omega^\zeta) = p_k^\zeta (\omega), \qquad k\in\mathbb{N}.
\end{equation}
To introduce an initial ensemble, take a sequence $(\omega^{(n)})_{n\in\mathbb{N}}$,  
$\omega^{(n)} = (\alpha^{(n)}, \beta^{(n)}, c^{(n)}) \in\Delta$, such that 
\begin{equation}\label{eq:4-45}
\alpha^{(n)}_{\zeta, 1} = O(1/\sqrt{n}), \quad \beta^{(n)}_{\zeta, 1} = O(1/\sqrt{n}) 
\quad (n\to\infty), \qquad \zeta\in\widehat{T}
\end{equation}
and set 
\begin{equation}\label{eq:4-46}
f_0^{(n)} = f_{\omega^{(n)}}\big|_{\mathfrak{S}_n(T)}, \qquad n\in\mathbb{N}.
\end{equation}
We assume that, along $\omega^{(n)}$ as $n\to\infty$, $\sqrt{n}$-rescaled 
Thoma mesures have weak convergence limit, namely, 
\begin{equation}\label{eq:4-47}
(\tau_{\omega^{(n)}}^\zeta)_{\sqrt{n}}(dx) = 
\tau_{\omega^{(n)}}^\zeta \bigl( \frac{1}{\sqrt{n}} dx) \ 
\xrightarrow[n\to\infty] \ \mathfrak{l}_0^\zeta(dx), \qquad \zeta\in\widehat{T}.
\end{equation}
Since $\mathrm{supp}(\tau_{\omega^{(n)}}^\zeta)_{\sqrt{n}}$ are uniformly 
bounded by \eqref{eq:4-45}, \eqref{eq:4-47} implies also moment convergence to 
compactly supported $\mathfrak{l}_0^\zeta$. 
We then realize the situation of Proposition\ref{prop:4-4} (and Corollary\ref{cor:4-3}, 
of course). 
In fact, let us see the assumptions in Theorem\ref{th:4-1}. 
\eqref{eq:4-46}, \eqref{eq:4-41}, \eqref{eq:4-44} and \eqref{eq:4-47} yield 
\begin{equation}\label{eq:4-49}
\lim_{n\to\infty} n^{\frac{k-1}{2}} f_0^{(n)} \bigl( \iota_{n,k} (k)_\theta \bigr) 
= \sum_{\zeta\in\widehat{T}} M_{k-1}(\mathfrak{l}_0^\zeta) \, 
\frac{\chi^\zeta_\theta}{(\dim\zeta)^k}, \qquad k\in\mathbb{N}.
\end{equation}
Then, \eqref{eq:4-6} immediately follows from multiplicativity of characters of 
$\mathfrak{S}_\infty (T)$. 
Comparing \eqref{eq:4-49} with \eqref{eq:4-5} and \eqref{eq:4-9}, we have 
\begin{equation}\label{eq:4-48}
R_{k+1}^\zeta = M_{k-1}(\mathfrak{l}_0^\zeta), \qquad k\in\mathbb{N}, \quad 
\zeta\in\widehat{T}
\end{equation}
In particular, $\lim_{n\to\infty} c^{(n)}_\zeta = R_2^\zeta = n^\zeta$. 
Since $(R_j^\zeta)$ determined by \eqref{eq:4-48} (with $R_1^\zeta =0$) is a 
conditionally positive-definite sequence, we obtain a $\boxplus$-infinitely 
divisible measure $\mathfrak{m}_0^\zeta$ related to L\'evy measure 
$\mathfrak{l}_0^\zeta$ through \eqref{eq:4-35} by following 
\cite[Theorem13.16, Proposition13.5]{NiSp06}. 
As further concrete specializations of \eqref{eq:4-45} and \eqref{eq:4-47}, 
we mention the following sequences of probabilities on $\mathbb{Y}_n(\widehat{T})$'s, 
or restriction of characters of $\mathfrak{S}_\infty(T)$ as \eqref{eq:4-46}, 
each corresponding to $\omega^{(n)} = (\alpha^{(n)}, \beta^{(n)}, c^{(n)})$ in 
\eqref{eq:4-40}. 
\begin{description}
\item[Plancherel] 
\[
\alpha^{(n)}_{\zeta, i} =0, \quad \beta^{(n)}_{\zeta, i} =0, \quad 
c^{(n)}_\zeta = \frac{(\dim\zeta)^2}{|T|}. \tag*{(P1)}
\]
The corresponding character is $\delta_e$. 

\item[Tensor type] \ Let 
$r, r^\prime, a_\zeta, b_\zeta, c_\zeta >0$ such that $a_\zeta+b_\zeta\leqq c_\zeta$, 
$\sum_\zeta c_\zeta =1$ be fixed, and 
$N, N^\prime\in\mathbb{N}$, \quad $c^{(n)}_\zeta = c_\zeta$, 
\begin{align*}
&\alpha^{(n)}_{\zeta, 1} = \cdots = \alpha^{(n)}_{\zeta, N} = \frac{a_\zeta}{N}, &  
&\alpha^{(n)}_{\zeta, N+1} = \alpha^{(n)}_{\zeta, N+2} = \cdots =0, & 
&N\sim r\sqrt{n}, \\ 
&\beta^{(n)}_{\zeta, 1} = \cdots = \beta^{(n)}_{\zeta, N^\prime} = 
\frac{b_\zeta}{N^\prime}, & 
&\beta^{(n)}_{\zeta, N^\prime +1} = \beta^{(n)}_{\zeta, N^\prime +2} = \cdots =0, & 
&N^\prime \sim r^\prime \sqrt{n}. \tag*{(P2)}
\end{align*}
If $T$ is trivial and $b=0$, $a=c=1$, the corresponding character is the one 
coming from the canonical action $U$ onto $(\mathbb{C}^M)^{\otimes n}$. 
In particular, the value of normalized trace of $U(k\text{-cycle})$ is $1/N^{k-1}$. 

\item[$q$-Plancherel] \ Let 
$r, r^\prime, a_\zeta, b_\zeta, c_\zeta >0$ such that $a_\zeta+b_\zeta\leqq c_\zeta$, 
$\sum_\zeta c_\zeta =1$ be fixed, and 
$0< q\leqq 1$, \ $0<q^\prime \leqq 1$, \quad $c^{(n)}_\zeta = c_\zeta$, 
\[ 
\alpha^{(n)}_{\zeta, i} = a_\zeta (1-q) q^{i-1}, \ 
\beta^{(n)}_{\zeta, i} = b_\zeta (1-q^\prime) q^{\prime\, i-1}, \quad 
1-q\sim \frac{1}{r\sqrt{n}}, \  1-q^\prime \sim \frac{1}{r^\prime\sqrt{n}}. \tag*{(P3)}
\] 
If $T$ is trivial and $b=0$, $a=c=1$, this is a $q$-analogue of the Plancherel measure 
in the sense that a symmetric group (algebra) is replaced by an associated Iwahori--Hecke algebra. 
\end{description}

\begin{table}[thb]
\begin{tabular}{ccccc} \toprule 
\  & (pausing time) & Plancherel & Tensor type & $q$-Plancherel \\ 
\midrule 
parameters $\omega^{(n)}$ & \ & (P1) & (P2) & (P3) \\ 
$R^\zeta(0,w)$ & \  & (R01) & (R02) & (R03) \\ 
$\mathfrak{l}^\zeta_0$ & \ & (L01) & (L02) & (L03) \\
$R^\zeta(t,w)$ \eqref{eq:4-28} & exponential-like & (RE1) & (RE2) & (RE3) \\
$\mathfrak{l}^\zeta_t$ \eqref{eq:4-36} & exponential-like & (LE1) & (LE2) & (LE3) \\
$R^\zeta(t,w)$ \eqref{eq:4-31}& $1/2$-stable & (RS1) & (RS2) & (RS3) \\ 
$\mathfrak{l}^\zeta_t$ \eqref{eq:4-37} & $1/2$-stable & (LS1) & (LS2) & (LS3) \\ \bottomrule 
\end{tabular}
\caption{$R^\zeta (t, \cdot)$ is Voiculescu's $R$-transform of probability 
$\mathfrak{m}^\zeta(t)$, and $\mathfrak{l}^\zeta_t$ is the associated L\'evy measure 
at $t\geqq 0$. 
For pausing time, \lq exponential-like\rq\ and \lq $1/2$-stable\rq\ correspond to 
Corollary\ref{cor:4-3} and Proposition\ref{prop:4-4} respectively.}
\label{tab:1}
\end{table}

\noindent 
In these three cases, we list their $R$-transforms $R^\zeta(0, w)$ (initial), $R^\zeta(t, w)$ 
(at time $t\geqq 0$), and L\'evy measures $\mathfrak{l}^\zeta_0$ (initial), 
$\mathfrak{l}^\zeta_t$ (at time $t\geqq 0$) in Table\ref{tab:1} filled out by the following formulas: 
\begin{align*}
&\frac{(\dim\zeta)^2}{|T|} w \tag*{(R01), (RE1), (RS1)} \\ 
&\frac{(\dim\zeta)^2}{|T|} \delta_0 \tag*{(L01), (LE1), (LS1)}
\end{align*}
\begin{align*}
&(c_\zeta- a_\zeta- b_\zeta) w + \frac{a_\zeta w}{1-(a_\zeta/r) w}+ 
\frac{b_\zeta w}{1+(b_\zeta /r^\prime) w} \tag*{(R02)} \\ 
&(c_\zeta- a_\zeta- b_\zeta) \delta_0 + a_\zeta \delta_{a_\zeta/r}+ 
b_\zeta \delta_{-b_\zeta/r^\prime} \tag*{(L02)}
\end{align*}
\begin{align*}
&\bigl\{ (1-e^{-t/m})\frac{(\dim\zeta)^2}{|T|} + (c_\zeta-a_\zeta-b_\zeta) e^{-t/m}
\bigr\} w + \frac{a_\zeta e^{-t/m} w}{1- (a_\zeta /r) e^{-t/m} w} + 
\frac{b_\zeta e^{-t/m} w}{1+ (b_\zeta /r^\prime) e^{-t/m} w}
\tag*{(RE2)} \\ 
&\bigl\{ (1-e^{-t/m})\frac{(\dim\zeta)^2}{|T|} + (c_\zeta-a_\zeta-b_\zeta) e^{-t/m}
\bigr\} \delta_0 \\ 
&\qquad\qquad\qquad\qquad\qquad + a_\zeta e^{-t/m} \delta_{(a_\zeta /r) e^{-t/m}} + 
b_\zeta e^{-t/m} \delta_{- (b_\zeta /r^\prime) e^{-t/m}}
\tag*{(LE2)}
\end{align*}
\begin{align*}
\Bigl[ \int_\mathbb{R} \frac{e^{-u^2/2}}{\sqrt{2\pi}} &\bigl\{ (1-e^{-\sqrt{t}|u|}) 
\frac{(\dim\zeta)^2}{|T|} + (c_\zeta-a_\zeta-b_\zeta) e^{-\sqrt{t}|u|} \bigr\} du\Bigr] w \\ 
&+ \int_\mathbb{R} \frac{e^{-u^2/2}}{\sqrt{2\pi}} \Bigl( 
\frac{a_\zeta e^{-\sqrt{t}|u|} w}{1- (a_\zeta /r) e^{-\sqrt{t}|u|} w} + 
\frac{b_\zeta e^{-\sqrt{t}|u|} w}{1+ (b_\zeta /r^\prime) e^{-\sqrt{t}|u|} w}\Bigr) du 
\tag*{(RS2)} \\ 
\Bigl[ \int_\mathbb{R} \frac{e^{-u^2/2}}{\sqrt{2\pi}} &\bigl\{ (1-e^{-\sqrt{t}|u|}) 
\frac{(\dim\zeta)^2}{|T|} + (c_\zeta-a_\zeta-b_\zeta) e^{-\sqrt{t}|u|} \bigr\} du \Bigr] 
\delta_0 (dx) \\ 
&+ \frac{\sqrt{2}}{\sqrt{\pi t}} \Bigl\{ r^\prime 
e^{-\frac{1}{2t}(\log \frac{-r^\prime x}{b_\zeta})^2} 1_{(-b_\zeta /r^\prime, 0)}(x) + 
r e^{-\frac{1}{2t}(\log \frac{rx}{a_\zeta})^2} 1_{(0, a_\zeta /r)}(x) \Bigr\} dx
\tag*{(LS2)}
\end{align*}
\begin{align*}
&(c_\zeta- a_\zeta- b_\zeta) w - r\log (1-\frac{a_\zeta w}{r})+ 
r^\prime \log (1+ \frac{b_\zeta w}{r^\prime}) 
\tag*{(R03)} \\ 
&(c_\zeta- a_\zeta- b_\zeta) \delta_0 (dx) + \bigl\{ r^\prime 1_{(-b_\zeta /r^\prime, 0)}(x) 
+ r 1_{(0, a_\zeta /r)}(x) \bigr\} dx 
\tag*{(L03)}
\end{align*}
\begin{align*}
&\bigl\{ (1-e^{-t/m})\frac{(\dim\zeta)^2}{|T|} + (c_\zeta-a_\zeta-b_\zeta) e^{-t/m}
\bigr\} w \\ 
&\qquad\qquad\qquad\qquad\qquad - r \log \bigl(1- \frac{a_\zeta e^{-t/m} w}{r} \bigr) + 
r^\prime \log \bigl(1+ \frac{b_\zeta e^{-t/m} w}{r^\prime} \bigr) 
\tag*{(RE3)} \\ 
&\bigl\{ (1-e^{-t/m})\frac{(\dim\zeta)^2}{|T|} + (c_\zeta-a_\zeta-b_\zeta) e^{-t/m}
\bigr\} \delta_0 (dx) \\ 
&\qquad\qquad\qquad\qquad\qquad + 
\bigl\{ r^\prime 1_{(-(b_\zeta /r^\prime) e^{-t/m}, 0)}(x) +
r 1_{(0, (a_\zeta /r) e^{-t/m})}(x) \bigr\}dx
\tag*{(LE3)} 
\end{align*}
\begin{align*}
&\Bigl[ \int_\mathbb{R} \frac{e^{-u^2/2}}{\sqrt{2\pi}} \bigl\{ (1-e^{-\sqrt{t}|u|}) 
\frac{(\dim\zeta)^2}{|T|} + (c_\zeta-a_\zeta-b_\zeta) e^{-\sqrt{t}|u|} \bigr\} du\Bigr] w \\ 
&\quad + \int_\mathbb{R} \frac{e^{-u^2/2}}{\sqrt{2\pi}} \bigl\{  
- r \log \bigl(1- \frac{a_\zeta e^{-\sqrt{t}|u|} w}{r} \bigr) + 
r^\prime \log \bigl(1+ \frac{b_\zeta e^{-\sqrt{t}|u|} w}{r^\prime} \bigr) \bigr\} du 
\tag*{(RS3)} \\ 
&\Bigl[ \int_\mathbb{R} \frac{e^{-u^2/2}}{\sqrt{2\pi}} \bigl\{ (1-e^{-\sqrt{t}|u|}) 
\frac{(\dim\zeta)^2}{|T|} + (c_\zeta-a_\zeta-b_\zeta) e^{-\sqrt{t}|u|} \bigr\} du\Bigr] 
\delta_0 (dx) \\ 
&\ + \Bigl\{ \Bigl( 
\int_{\frac{1}{\sqrt{t}}\log (-\frac{r^\prime x}{b_\zeta})}^{\frac{1}{\sqrt{t}}
\log (-\frac{b_\zeta}{r^\prime x})}
\frac{e^{-u^2/2}}{\sqrt{2\pi}} du \Bigr) r^\prime\, 1_{(-b_\zeta /r^\prime, 0)}(x) + 
\Bigl( 
\int_{\frac{1}{\sqrt{t}}\log \frac{rx}{a_\zeta}}^{\frac{1}{\sqrt{t}}\log \frac{a_\zeta}{rx}}
\frac{e^{-u^2/2}}{\sqrt{2\pi}} du \Bigr) r\, 1_{(0, a_\zeta /r)}(x) \Bigr\} dx.
\tag*{(LS3)}
\end{align*}
We note that consistency of L\'evy measures at $t=0$ can be understood through weak convergence of 
$\mathfrak{l}^\zeta_t$ to $\mathfrak{l}^\zeta_0$ as $t\to 0$.

\bigskip

\textbf{Acknowledgments} \ 
The author expresses deep appreciation to Professor Takeshi Hirai 
for long-term support, instruction, and collaboration in our research projects.

\end{document}